\documentclass[12pt,reqno]{amsart}
\usepackage{amssymb,latexsym}
\usepackage{graphicx}
\usepackage{url}		%does nice formatting of URLs

\calclayout

\makeatletter
\newcommand{\monthyear}[1]{%
  \def\@monthyear{\uppercase{#1}}}
\newcommand{\volnumber}[1]{%
  \def\@volnumber{\uppercase{#1}}}
\AtBeginDocument{%
\def\ps@plain{\ps@empty
  \def\@oddfoot{\@monthyear \hfil \thepage}%
  \def\@evenfoot{\thepage \hfil \@volnumber}}
\def\ps@firstpage{\ps@plain}
\def\ps@headings{\ps@empty
  \def\@evenhead{%
    \setTrue{runhead}%
    \def\thanks{\protect\thanks@warning}%
    \uppercase{}\hfil}%
  \def\@oddhead{%
    \setTrue{runhead}%
    \def\thanks{\protect\thanks@warning}%
    \hfill\uppercase{On Some Primality Tests based on Binomial}}%
  \let\@mkboth\markboth
  \def\@evenfoot{%
    \thepage \hfil \@volnumber}%
  \def\@oddfoot{%
    \@monthyear \hfil \thepage}%
  }%
\footskip=25pt
\pagestyle{headings}%
}
\makeatother

\theoremstyle{plain}
\numberwithin{equation}{section}
\newtheorem{thm}{Theorem}[section]
\newtheorem{theorem}[thm]{Theorem}
\newtheorem*{lucthm}{Lucas' Theorem} % not numbered
\theoremstyle{conjecture}

\theoremstyle{corollary}

\theoremstyle{remark} % 'style changed again'
\newtheorem{remark}[thm]{Remark}
\newtheorem{proposition}[thm]{Proposition}

\newcommand{\bq}{\begin{equation}}
\newcommand{\eq}{\end{equation}}
\newcommand{\bareq}{\begin{eqnarray}}
\newcommand{\eareq}{\end{eqnarray}}

\begin{document}
\vspace*{-0.5 in}

%% replace the values in the next three lines by the correct information
\monthyear{Month Year}
\volnumber{Volume, Number}
\setcounter{page}{1}

\title{On Primality Tests Grounded on Binomial Coefficients}

\author{Dario T. de Castro
}
\address{Instituto Federal do Rio de Janeiro - Campus Nil\'opolis \\ Rua Cel. D\'elio Menezes Tavares 1045 - Nil\'opolis \\ Rio de Janeiro - RJ - Brazil - CEP 26530-060
}
\email{dario.neto@ifrj.edu.br
}
\thanks{}

\begin{abstract}
In this paper, we introduce two primality tests based on new divisibility properties of binomial coefficients. These new properties were enunciated and proved in previous work. We also study two similar tests that can be obtained from well-known results in Number Theory. At the end we compare our results with the existing ones.
\end{abstract}

\maketitle

\section{Introduction}
In a previous work \cite{Castro} we introduced some results related to the divisibility properties of binomial coefficients. In particular, we proved a result that is equivalent to the following proposition:
\begin{proposition}\label{prop} Given a prime number $p$ and positive integers $n$ and $j$, we have
\bq \label{eq1}
{n - 1 \choose p^j - 1} \equiv \begin{cases}
 1 \;\mbox{\rm (mod }p\mbox{\rm )} & \text{if } \; p^j \vert n \\
 0 \;\mbox{\rm (mod }p\mbox{\rm )}  & \text{if } \; p^j \nmid n. \\
\end{cases}
\eq
\end{proposition}
From this congruence relation we defined a Boolean operator related to whether or not $p ^ j$ divides $n$. We have also shown that the primality of $p$, as it appears in (\ref{eq1}), is a necessary condition if this formula is expected to hold for all $n \in \mathbb{N^*}$. In fact, we shall prove in what follows that, for a composite number $q > 0$, divisible by a given prime number $\hat p$, ${q + \hat p - 1 \choose q - 1} \not\equiv 0 \;\mbox{\rm (mod }q\mbox{\rm )}$ even though $q \nmid q + \hat p$.

\begin{theorem}\label{thm1} Let $q > 0$ be a composite number divisible by the prime $\hat p$. If $n_{q} = q + \hat p$, then ${n_q - 1 \choose q - 1} \not\equiv 0 \;\mbox{\rm (mod }q\mbox{\rm )}$ even though we have $q \nmid n_{q}$.\end{theorem}
\begin{proof} Initially, we have:
\begin{align}
{q + \hat p - 1  \choose q - 1} &= {q + \hat p - 1  \choose \hat p} \nonumber \\
&= {q + \hat p - 1  \choose \hat p - 1} \frac{q}{\hat p}.
\end{align}
Since ${q + \hat p - 1  \choose \hat p - 1} \equiv 1 \; \mbox{\rm (mod }\hat p\mbox{\rm )}$, as can be deduced from either Lucas' theorem (see Appendix) or Proposition \ref{prop}, we have
\bq
{q + \hat p - 1  \choose q - 1} = \left(M \hat p + 1 \right) \frac{q}{\hat p},
\eq
where $M$ is a positive integer. Thus, given that $q \nmid q/\hat p$, we conclude that
\bq
{q + \hat p - 1 \choose q - 1} \not\equiv 0 \;\mbox{\rm (mod }q\mbox{\rm )}.
\eq
\end{proof}
\begin{remark} If a prime $p$ does not divide the composite $q$, then ${q + p - 1  \choose q - 1} \equiv 0 \; \mbox{\rm (mod } q\mbox{\rm )}$. We leave it to the reader to verify this fact. \end{remark}

\begin{remark} Theorem \ref{thm1} is a particular case of a more general statement proved in our previous work \cite{Castro}. \end{remark}

The following theorem is a well known result in number theory and serves as basis for some tests of primality \cite{AKS}. We state it here without proof (see \cite{Ogilvy}).
\begin{theorem}\label{thm2} A positive integer $n$ is a prime number if and only if it divides ${n \choose k}$, for all $k$ satisfying $1 \le k \le n-1$.
\end{theorem}
From this property one can write a primality test on any positive integer, which is given in terms of binomial coefficients. It can be stated as follows:

\noindent
{\it A positive integer $n$ is prime if and only if}
\begin{equation}\label{pt0a}
\sum_{i=1}^{n - 1}{n \choose i} \;\mbox{\rm mod }n = 0.
\end{equation}
However, as shown in \cite{Ogilvy}, if $n$ is composite and $p$ is the smallest prime factor of $n$, then
\begin{equation}
{n \choose p} \;\mbox{\rm mod }n > 0.
\end{equation}
Thus, one can improve the test above, so that it now reads as

\noindent
{\it A positive integer $n > 3$ is prime if and only if}
\begin{equation}\label{pt0a1}
\sum_{i=1}^{\pi(\sqrt{n})}{n \choose p_i} \;\mbox{\rm mod }n = 0,
\end{equation}
where $p_i$ is the $i$th prime number and $\pi(x)$ represents the prime counting function, which yields the number of primes that are less than or equal to $x$.

In connection with Theorem \ref{thm2} and as a direct application of the Pascal's rule for binomial coefficients \cite{Cannon}, i.e.,
\begin{equation}
{n \choose k} = {n - 1 \choose k} + {n - 1 \choose k - 1},
\end{equation}
an alternative primality test can be enunciated, also in terms of binomial coefficients, as follows:

\noindent
{\it A positive integer $n > 1$ is prime if and only if}
\begin{equation}\label{pt0b}
\sum_{i=0}^{\left \lfloor \frac{n}{2}-1 \right \rfloor}{n + i\choose i + 1} \;\mbox{\rm mod }n = 0,
\end{equation}
where, for a real number x, we denote by $\lfloor x \rfloor$ the largest integer less than or equal to x.

In section 2 we shall present and prove two primality tests based on Proposition \ref{prop} and Theorem \ref{thm1}, which represent some improvement over (1.7) and (1.9).

\section{Main Results}
Our first primality test follows immediately from Proposition \ref{prop}.
\begin{theorem}\label{thm21} Let $n > 3$ be an integer and $p_i$ be the $i$th prime number. Let $\pi(x)$ be the prime counting function. Then, $n$ is prime if and only if
\begin{equation}\label{test1}
\sum_{i=1}^{\pi(\sqrt{n})}{n - 1 \choose p_i - 1} \;\mbox{\rm mod }p_i = 0
\end{equation}
\end{theorem}
\begin{proof}
It is clear that each term of the sum in (\ref{test1}) is an application of Proposition \ref{prop}, with $j = 1$, for a given prime which is less than or equal to $\sqrt{n}$.
\end{proof}

Now, we shall introduce a second primality test, which follows from Proposition \ref{prop} and Theorem \ref{thm1}.
\begin{theorem}\label{thm22} Let $n > 3$ be an integer and $p_i$ be the $i$th prime number. Let $\pi(x)$ be the prime counting function. Then, $n$ is prime if and only if
\begin{equation}\label{test2}
\sum_{i=1}^{\pi(\sqrt{n})}{n + p_i - 1 \choose n - 1} \;\mbox{\rm mod }n = 0
\end{equation}
\end{theorem}
\begin{proof} Since $p_i < n$, it follows that $n \nmid n + p_i$. If $n$ is prime, all terms in the sum will be zero as implied by Proposition \ref{prop}. On the other hand, if $n$ is composite, at least one of its prime factors must be in the set ${\mathcal P}_n = \{p_i \; \vert \; 1 \leq i \leq \pi(\sqrt{n})\}$. Therefore, according to Theorem \ref{thm1}, whenever $p_i \in {\mathcal P}_n$ is a prime factor of $n$, a positive term will be generated and the sum will be greater than zero. \end{proof}

\begin{remark} We believe the tests presented in Theorems \ref{thm21} and \ref{thm22} are of theoretical interest since: (a) they are fundamentally combinatorial as they are based exclusively on the divisibility properties of binomial coefficients, (b) they are deterministic, (c) for $n$ integer and greater than $3$, they are unconditionally correct, (d) one can immediately get from them the list of prime divisors of $n$ which are less than or equal to $\sqrt{n}$, and (e) when compared to (\ref{pt0a1}) and (\ref{pt0b}), they generally allow us to both deal with smaller numbers and perform fewer operations to achieve a result.\end{remark}

\begin{remark} Different algorithms can be used to evaluate the sums of remainders in (\ref{test1}) and (\ref{test2}), each of which presenting a different level of computational complexity. Although we do not provide in this work a detailed analysis of the efficiency of these two tests, preliminary estimates point to asymptotic time complexities which depend on the exponential of the number of bits required to represent $n$ in binary \cite{Granville}.\end{remark}

\section{Appendix}
\begin{lucthm} Let $p$ be a prime number and let $R$ and $S$ be positive integers such that
\bq
R = r_0 + r_1 p + r_2 p^2 + ... + r_m p^m \, , \; \; r_i \in \{0,1,2,...,p-1\},
\eq
and
\bq
S = s_0 + s_1 p + s_2 p^2 + ... + s_l p^l \, , \; \; s_i \in \{0,1,2,...,p-1\},
\eq
where $r_i$ and $s_i$ are, respectively, the digits of $R$ and $S$ when written in base $p$. We then have
\bq \label{Lucas}
{R \choose S} \equiv  \prod_{i = 0}^{\max\{m,l\}}\, {r_i \choose s_i} \; \mbox{\rm (mod }p\mbox{\rm )}.
\eq
We adopt here the convention ${r \choose s} = 0$ if $s$ is either greater than $r$ or smaller than zero.
\end{lucthm}

\section{Acknowledgment}
The author wishes to thank Professors Diego Marques and Reinaldo de Melo for useful and inspiring suggestions.

\end{document}